\documentclass[conference]{IEEEtran}
\IEEEoverridecommandlockouts
\usepackage{cite}
\usepackage{amsmath,amssymb,amsfonts}
\usepackage{algorithmic}
\usepackage{graphicx}
\usepackage{graphics}
\usepackage{textcomp}
\usepackage{xcolor}
\usepackage{comment}
\usepackage{caption}
\usepackage{epstopdf}
\usepackage{subcaption}
\usepackage{tablefootnote}
\usepackage[hidelinks]{hyperref}

\usepackage[linesnumbered,ruled,vlined]{algorithm2e}
\SetKwInput{KwInput}{Input}                
\SetKwInput{KwOutput}{Output} 
\SetKwInput{KwInit}{Initialize}  
\SetKw{Break}{break}

\def\BibTeX{{\rm B\kern-.05em{\sc i\kern-.025em b}\kern-.08em
    T\kern-.1667em\lower.7ex\hbox{E}\kern-.125emX}}
\begin{document}

\title{Manifold Quadratic Penalty Alternating Minimization for Sparse Principal Component Analysis}

\author{\IEEEauthorblockN{Tarmizi Adam}
\IEEEauthorblockA{\textit{ViCubeLab, Faculty of Computing} \\
\textit{Universiti Teknologi Malaysia}\\
Jalan Iman, 81310, Skudai, Johor, Malaysia \\
tarmizi.adam@utm.my}

}

\maketitle

\begin{abstract}
Optimization on the Stiefel manifold or with orthogonality constraints is an important problem in many signal processing and data analysis applications such as Sparse Principal Component Analysis (SPCA). Algorithms such as the Riemannian proximal gradient algorithms addressing this problem usually involve an intricate subproblem requiring a semi-smooth Newton method hence, simple and effective operator splitting methods extended to the manifold setting such as the Alternating Direction Method of Multipliers (ADMM) have been proposed. However, another simple operator-splitting method, the Quadratic Penalty Alternating Minimization (QPAM) method which has been successful in image processing to our knowledge, has not yet been extended to the manifold setting. In this paper, we propose a manifold QPAM (MQPAM) which is very simple to implement. The iterative scheme of the MQPAM consists of a Riemannian Gradient Descent (RGD) subproblem and a subproblem in the form of a proximal operator which has a closed-form solution. Experiments on the SPCA problem show that our proposed MQPAM is at par with or better than several other algorithms in terms of sparsity and CPU time.
\end{abstract}


\section{Introduction}
This paper is concerned with the nonsmooth nonconvex  Riemannian optimization of the form: 
\begin{equation}
\begin{aligned}
&\underset{\mathbf{X}}{minimize}\, f\left( \mathbf{X}\right) + g\left(\mathbf{X} \right) \\
& \text{\textit{subject to}}\, \mathbf{X} \in \mathcal{M}, \\
\end{aligned}
\label{eq: probMain}
\end{equation}where $\mathcal{M} = \text{St}\left(n,\,p \right)=\{ \mathbf{X} \in \mathbb{R}^{n \times p}: \mathbf{X}^\top \mathbf{X} = \mathbf{I}_p \}$ is the Stiefel manifold \cite{boumal2023introduction,absil2008optimization}. The function $f$ is smooth and possibly nonconvex, while $g$ is nonsmooth and convex.

Problems of the form (\ref{eq: probMain}) arise in many applications such as unsupervised feature selection \cite{yang2011,tang2012}, dictionary learning \cite{sun2016b,sun2016a},   sparse blind deconvolution \cite{zhang2017global},  compressed modes \cite{ozolicnvs2013}, and sparse principal component analysis \cite{jolliffe2003,huang2022,xiao2021exact}.

An important data analysis tool is the principal component analysis (PCA) \cite{pearson1901,hotelling1933}. The PCA is used to find a low-dimensional representation of a high-dimensional dataset. However, the loading vectors obtained in PCA are typically non-sparse and would result in difficult data interpretability in certain applications requiring sparse loading vectors. To enforce sparsity in the loading vectors, Sparse Principal Component Analysis (SPCA) can be modeled as,
\begin{equation}
\begin{aligned}
&\underset{\mathbf{X}}{minimize}\, F\left( \mathbf{X} \right)= - \frac{1}{2}\text{Tr}\left( \mathbf{X}^\top \mathbf{A}^\top \mathbf{A}\mathbf{X} \right) +  \mu\| \mathbf{X} \|_1 \\
&subject\,to\quad \mathbf{X}\in \text{St}\left(n,\,p \right),
\end{aligned}
\label{eq:spca1}
\end{equation}
where $\text{Tr}\left( \cdot\right)$ is the trace operator of a matrix, $\mathbf{A}\in\mathbb{R}^{m \times n}$ is a data matrix, and $\mu$ is the regularization parameter. Model (\ref{eq:spca1}) \cite{jolliffe2003} imposes sparsity and orthogonality on the loading vectors as a result of the $\ell_1$-norm regularization and the Stiefel manifold constraint. Obviously, (\ref{eq:spca1}) is in the form of (\ref{eq: probMain}).

In general, manifold-constrained optimization has been a hot topic for the past decades due to its wide applicability in applications such as SPCA mentioned above \cite{hu2020brief,boumal2023introduction,huang2022}. Furthermore, algorithmic developments for manifold optimization such as the Riemannian Gradient descent (RGD) \cite{boumal2023introduction}, Riemannian Conjugate Gradient method (RCG) \cite{sato2021riemannian}, Riemannian Proximal Gradient (RPG) \cite{huang2022riemannian,chen2020proximal,huang2021robust}, Manifold Alternating Direction Method of Multipliers (MADMM) \cite{kovnatsky2016madmm}, and Riemannian ADMM (RADMM) \cite{li2022riemannian} has also contributed in increasing the interest in the field.

In this paper, we develop a manifold operator splitting method called the Manifold Quadratic Penalty Alternating Minimization (MQPAM) for problems of the form (\ref{eq: probMain}) and focus its application on solving the SPCA problem (\ref{eq:spca1}). Different from other operator splitting schemes such as the MADMM and RADMM \cite{kovnatsky2016madmm,li2022riemannian}, our method is simple to implement and does not involve optimization subproblems in the dual space (Lagrange multiplier updates).

The MQPAM is motivated and the extension of the quadratic penalty alternating minimization (QPAM) method which is used successfully in image processing and sparse signal recovery \cite{wang2008new,xie2018new,xu2011image}.

When comparing with Riemannian proximal gradient methods such as \cite{chen2020proximal,huang2022} the MQPAM does not rely on complex inner iterations such as the semismooth Newton method but only relies on an RGD \cite{boumal2023introduction} and a closed-form proximal step. 
This results in a very fast algorithm for SPCA with comparable or even better results to the aforementioned methods.

The rest of the paper is structured as follows. First, we overview works related to our proposed method in Section \ref{sec:RW}. Necessary background is presented in Section \ref{sec:Pre}. The proposed algorithm is given in Section \ref{sec:PM}. Finally, results and conclusions are presented in Sections \ref{sec:ER} and \ref{sec:C}.

\section{Related Works}
\label{sec:RW}
Solving optimization problems with orthogonality constraints or on the Stiefel manifold has gained interest in practical applications. Many algorithms for orthogonality-constrained problems have been proposed \cite{lai2014splitting,siegel2021accelerated,chen2020proximal,chen2016augmented,li2022riemannian}. Among these algorithms, specific to problems (\ref{eq: probMain}) and (\ref{eq:spca1}) are the class of operator splitting methods and proximal gradient methods\cite{kovnatsky2016madmm,li2022riemannian,lai2014splitting,chen2020proximal,huang2022}. Originally, these algorithms were studied in the Euclidean space and were extended to the manifold setting \cite{boyd2011distributed,goldstein2009split,beck2009fast}.

The manifold proximal gradient method and its accelerated version for the SPCA model (\ref{eq:spca1}) have been recently proposed \cite{chen2020proximal,huang2022}. However, a shortcoming of these algorithms is the complex inner iteration they possess i.e., they require a semi-smooth Newton inner iteration hence, increasing the computational time. This has been shown in the results of \cite{deng2023manifold} when compared to Augmented Lagrangian (AL) and operator-splitting methods such as the ADMM \cite{li2022riemannian}. 

The manifold MADMM has been studied and proposed in \cite{kovnatsky2016madmm} and quite recently, the RADMM has been proposed in \cite{li2022riemannian}. Both MADMM and RADMM have shown good performance in terms of CPU time for several data analysis and computer vision applications.

The QPAM algorithm is another operator-splitting method that has shown success in image and signal processing \cite{wang2008new,xie2018new,tan2014smoothing}. Unlike ADMM-type methods, QPAM does not construct the augmented Lagrangian function of the original problem. Instead, it adds a simple quadratic penalty function. Consequently, QPAM iterations are much easier compared to ADMM i.e., it does not need a dual (Lagrange multipliers) update step. Due to its simplicity and the success of operator splitting methods such as the MADMM and the RADMM, in this paper, we are motivated to extend the QPAM to the manifold setting. To this end, the contributions of this paper are:
\begin{enumerate}
\item We propose a manifold QPAM (MQPAM) for optimization problem (\ref{eq: probMain}) hence, extending the QPAM to the Riemannian setting.
\item We show that the proposed MQPAM is simple to implement consisting of a Riemannian Gradient Descent (RGD) step and a proximal step. Furthermore, we show that the proximal step of the MQPAM is closely linked to the Moreau smoothing technique.
\item We conduct experiments and discussions on the SPCA problem (\ref{eq:spca1}) and show that the MQPAM gives better or comparable results compared to other competing methods on the same problem.
\end{enumerate}

\section{Preliminaries}
\label{sec:Pre}

This paper relies on several notations. The Euclidean inner product between two matrices $\mathbf{A}$ and $\mathbf{B}$ is denoted by $\langle \mathbf{A},\,\mathbf{B} \rangle = \text{Tr}\left(\mathbf{A}^\top \mathbf{B} \right)$ where $\text{Tr}\left(\cdot \right)$ denotes the trace operator. The $\ell_1$ and $\ell_2$ norms are denoted as $\| \cdot \|_1$, and $\| \cdot\|_2$ respectively.

The tangent space to a manifold $\mathcal{M}$ at point $\mathbf{X}$ is noted as $T_\mathbf{X} \mathcal{M}$. For the Stiefel manifold $\text{St}\left(n,\,p \right)$ its tangent space is given by $T_\mathbf{X} \mathcal{M} = \{\mathbf{V} \in \mathbb{R}^{n \times p}:\mathbf{X}^\top\mathbf{V} + \mathbf{V}^\top \mathbf{X} = 0   \}$. The Riemannian metric on $T_\mathbf{X} \mathcal{M}$ induced from the Euclidean inner product is $\langle \zeta^\top ,\, \eta \rangle = \text{Tr}\left( \zeta^\top ,\, \eta \right)$ such that  $\forall \zeta,\, \eta \in T_\mathbf{X} \mathcal{M}$. Other notations will be defined when they occur throughout the paper.

\section{Proposed Method and Algorithm}
\label{sec:PM}

We start by proposing our MQPAM for solving (\ref{eq: probMain}). To split the smooth and nonsmooth part, we introduce a new variable $\mathbf{Y}$ and add a quadratic penalty term to (\ref{eq: probMain}):
\begin{equation}
\begin{aligned}
&\underset{\mathbf{X},\, \mathbf{Y}}{minimize}\, f\left( \mathbf{X}\right) + g\left(\mathbf{Y} \right) + \frac{\beta}{2}\| \mathbf{Y} - \mathbf{X} \|_2^2 \\
& \text{\textit{subject to}}\, \mathbf{Y} = \mathbf{X},\, \mathbf{X} \in \mathcal{M}. \\
\end{aligned}
\label{eq: probSplit1}
\end{equation}
This simple and computationally fast quadratic penalty splitting traces its root back to Courant \cite{courant1943} and has been proposed in \cite{wang2008new} in the context of Total Variation image restoration.

To solve (\ref{eq: probSplit1}), the MQPAM then iteratively solves,
\begin{equation}
\begin{cases}
\mathbf{X}_{k+1} = \underset{\mathbf{X} \in \mathcal{M}}{arg\,min}\, f\left( \mathbf{X}\right) + \frac{\beta}{2}\| \mathbf{Y}_k - \mathbf{X} \|_2^2 \\
\mathbf{Y}_{k+1} = \underset{\mathbf{Y}}{arg\,min}\, g\left(\mathbf{Y} \right) + \frac{\beta}{2}\| \mathbf{Y} - \mathbf{X}_{k+1} \|_2^2. \\
\end{cases}
\label{eq:ramItr}
\end{equation}

From the iterative scheme (\ref{eq:ramItr}), the $\mathbf{X}$ subproblem is a smooth and differentiable Riemannian optimization problem constrained on the Stiefel manifold. For this, the RGD can be utilized to solve it as follows \cite{boumal2023introduction},
\begin{equation}
\mathbf{X}_{k+1} = \text{Retr}_\mathbf{X}  \left( -\eta \mathcal{G}_k \right),
\label{eq:rgd}
\end{equation}
where $\mathcal{G}_k = \text{Proj}_{T_\mathbf{X} \mathcal{M}}\left(\mathbf{G}_k\right)$ is the Riemannian gradient at the $k^{\text{th}}$ iteration, $\mathbf{G}_k = \nabla f\left( \mathbf{X}\right) + \beta \left( \mathbf{X} - \mathbf{Y}_k \right)$ is the Euclidean gradient of the left-hand side of the $\mathbf{X}$ subproblem in (\ref{eq:ramItr}), and $\eta$ is the RGD step size.
Here, projection $\text{Proj}_{T_\mathbf{X}  \mathcal{M}} \left(\cdot \right)$ defines the projection onto the tangent space $T_\mathbf{X}  \mathcal{M}$ of $\text{St}\left(n,\,p \right)$ defined as \cite{huang2022}:
\[
\text{Proj}_{T_\mathbf{X}  \mathcal{M}} \left(\mathbf{G} \right) = \mathbf{G} - \mathbf{X} \text{skew}\left(\mathbf{X}^\top \mathbf{G} \right),
\]
where $\text{skew}\left(\mathbf{X}^\top \mathbf{G}\right)=\frac{ \mathbf{X}^\top \mathbf{G} - \mathbf{G}^\top \mathbf{X}     }{2}$ is the skew-symmetric part of the matrix \cite{boumal2023introduction}.
The retraction $\text{Retr}_\mathbf{X}  \left(\cdot\right)$, from the tangent space $T_\mathbf{X}  \mathcal{M}$ to the manifold $\mathcal{M}$ can be done in several ways. This includes QR decomposition, polar decomposition, and the Cayley transform \cite{boumal2023introduction,absil2008optimization}. In this work, we choose the polar decomposition given by,
\[
\text{Retr}_\mathbf{X}  \left(\mathcal{G}\right) = \left( \mathbf{X} + \mathcal{G} \right) \left( \mathbf{I}_p + \mathcal{G}^\top \mathcal{G} \right)^{- \frac{1}{2} }.
\]

The $\mathbf{Y}$ subproblem can be solved via standard proximal operators of the function $g\left( \cdot \right)$ defined as:
\begin{equation}
\text{prox}_{\mu g}\left( \mathbf{U} \right) = \underset{\mathbf{Z}}{\operatorname{arg\,min}}\,\, \frac{1}{2\mu}\| \mathbf{Z} - \mathbf{U} \|^2_2 + g\left(\mathbf{Z} \right).
\label{eq:prox_lam}
\end{equation}
As an example, if $g\left( \cdot \right)$ is the $\ell_1$-norm, then this subproblem can be solved via the soft-thresholding operator \cite{wright2022optimization}. 

To show this, let us use a smoothing strategy using the Moreau envelope of the function $g\left(\cdot\right)$ with parameter $\mu = \frac{1}{\beta}$,

\begin{equation}
g_\mu\left( \mathbf{U}\right) = \underset{\mathbf{Z}}{ \inf} \, \frac{1}{2\mu}\| \mathbf{Z} - \mathbf{U} \| ^2_2 + g\left(\mathbf{Z} \right).
\label{eq:morEnv}
\end{equation}
Now, the $\mathbf{Y}$ subproblem in (\ref{eq:ramItr}) can be seen as:

\begin{equation}
g_{\mu}\left(\mathbf{X}_{k+1} \right) = \underset{\mathbf{Y}}{\inf}\, \frac{1}{2\mu} \|\mathbf{Y} - \mathbf{X}_{k+1} \|_2^2 + g\left( \mathbf{Y} \right),
\label{eq:morEnv}
\end{equation}
i.e., minimizing a smoothed diffirentiable version of $g\left(\cdot\right)$ with Lipschitz constant $L=\beta$ and is closely related to the proximal operator (\ref{eq:prox_lam}). Furthermore, it is known that the proximal operator can be interpreted as minimizing a proximal linearization term. In the case of (\ref{eq:morEnv}), let $h\left( \mathbf{Y} \right) = \frac{1}{2\mu} \| \mathbf{Y} - \mathbf{X}_{k+1} \|_2^2$ then the $\mathbf{Y}$ subproblem in terms of minimizing the proximal linearization is,

\begin{equation}
\begin{aligned}
    \mathbf{Y}_{k+1} = \underset{\mathbf{Y}}{arg\,min}\, h\left( \mathbf{X}_{k+1}\right) &+ \langle \nabla h\left( \mathbf{X}_{k+1} \right),\, \mathbf{Y} - \mathbf{X}_{k+1} \rangle \nonumber \\
    &+ h\left( \mathbf{Y} \right) + g\left( \mathbf{Y} \right). \nonumber \\
\end{aligned}
\end{equation}
This implies,
 \begin{equation}
 	\mathbf{Y}_{k+1} = \underset{\mathbf{Y}}{arg\,min}\, \frac{1}{2\mu} \|\mathbf{Y} - \left( \mathbf{X}_{k+1} - \frac{1}{\mu} \nabla h\left( \mathbf{X}_{k+1}\right) \right)\|_2^2 + g\left( \mathbf{Y} \right),
 \end{equation}
which by the definition of the proximal operator (\ref{eq:prox_lam}) gives,
\begin{equation}
\mathbf{Y}_{k+1} = \text{prox}_{\mu g} \left( \mathbf{X}_{k+1} - \frac{1}{\mu} \nabla h\left( \mathbf{X}_{k+1}\right) \right),
\end{equation}
implying that the quadratic penalty method \cite{wang2008new} can be interpreted as a Moreau smoothing strategy.
This smoothing strategy is used in several previous works to obtain fast and simple computable subproblems \cite{tan2014smoothing,xie2018new,li2022riemannian}. As a result, the proposed MQPAM also inherits a simple and fast computation of the $\mathbf{Y}$ subproblem.
The overall MQPAM algorithm for solving (\ref{eq: probSplit1}) is listed in Algorithm \ref{al:rqpam_1}.

\begin{algorithm}
\caption{Manifold Quadratic Penalty Alternating Minimization (MQPAM)}
\label{al:rqpam_1}

\KwInit{$\mathbf{X}_0$, $\mathbf{Y}_0$, $\eta > 0$, $tol$, $Nit$, $Rgdit$}

\For{$k=0$ \KwTo $Nit$}
{
    \tcp{Start RGD iterations}
\For{$j=0$ \KwTo $Rgdit$}   
{
$\mathbf{G}_{k+1} = \nabla f\left( \mathbf{X}_k\right) + \beta \left( \mathbf{X}_k - \mathbf{Y}_k \right)$ \\
$\mathcal{G}_{k+1} = \text{Proj}_{T_X \mathcal{M}}\left(\mathbf{G}_k\right)$ \\
\If{$\| \mathcal{G}_{k+1} \|_2 < tol$}
{
Break
}
$\mathbf{X}_{k+1} = \text{Retr}_X \left( -\eta \mathcal{G}_{k+1} \right)$ \\
}

$\mathbf{Y}_{k+1} = \text{prox}_{\mu g} \left( \mathbf{X}_{k+1} - \frac{1}{\mu} \nabla h\left( \mathbf{X}_{k+1}\right) \right)$

$k= k+1$
}
\end{algorithm}

Several remarks ensue in Algorithm 1. First, the MQPAM is very simple to implement with one RGD call and one proximal operation which usually involves closed-form solutions. For the RGD, it is shown that the per-iteration complexity of the RGD obtaining an $\epsilon$-stationary point is of order $\mathcal{O}\left( \frac{1}{\epsilon^2} \right)$ \cite{boumal2019,bento17}.
Secondly, as mentioned, we have chosen to compute the retraction using the polar decomposition which costs $3np^2+\mathcal{O}\left( p^3\right)$ flops \cite{chen2020proximal}. For this, step 7 of Algorithm \ref{al:rqpam_1} involves computing the Singular Value Decomposition (SVD) of the matrix $\mathbf{X} + \mathcal{G}$ on each iteration $k$ \cite{boumal2023introduction}.

\begin{table*}[t]
	\caption{Comparative analysis between several splitting methods. Due to spacing and different $\mu$ values for MQPAM, the results for MQPAM are separated in Table \ref{tab:resMam}. }
	\centering
    \resizebox{18cm}{!}{
	\begin{tabular}{ccccccccccccccccc}
		\hline
		&  &  &  &SOC & &  &&  & MADMM  & & & & & RADMM  \\ 
		\hline
		&$\mu$  & $n/p$  &Obj  &Err  &CPU  &Spar &  &Obj  &Err  &CPU  &Spar & &Obj  &Err  &CPU  &Spar   \\
		\hline
		&  &$300/50$  &$31.361$  &$0.0018$   &$0.0958$ &$0.9966$ &  &$31.519$ &$0.00303$ &$0.0703$ &$0.9606$ & &$31.6981$ &$0.00132$ &$0.0585$ &$0.9572$  \\
		&$1$  &$300/100$  &$61.515$ &$0.0018$3  &$0.1023$ &$0.9966$ &  &$61.604$ &$0.00091$ &$0.2457$ &$0.9565$ & &$62.7107$ &$0.00179$ &$0.1905$ &$0.9635$  \\
		&  &$500/50$  &$39.187$  &$0.00205$  &$0.2662$  &$0.9979$  & &$39.014$  &$0.00102$ &$0.1728$ &$0.9689$  & &$38.8266$ &$0.00106$ &$0.1685$ &$0.9574$   \\
		&  &$500/100$  &$76.805$  &$0.00141$  &$0.3631$  &$0.9979$  & &$77.101$  &$0.00058$ &$0.4003$ &$0.9657$ & &$75.6348$ &$0.00109$ &$0.3788$ &$0.9575$  \\
		\hline
		&  &\textbf{Average}  &\textbf{52.217}  &\textbf{0.00177}  &\textbf{0.2069} &\textbf{0.9973}  &  &\textbf{52.309} &\textbf{0.00139} &\textbf{0.2223} &\textbf{0.9629} & &\textbf{52.218} &\textbf{0.00132} &\textbf{0.1991} &\textbf{0.9589}  \\
		\hline\hline
		&  &$300/50$  &$22.662$  &$0.00222$   &$0.0794$ &$0.9965$ &  &$19.154$ &$0.00079$  &$0.1155$ &$0.9538$ & &$22.3477$ &$0.00214$ &$0.0598$ &$0.9469$ \\
		&$0.8$  &$300/100$  &$42.763$  &$0.00141$   &$0.1442$ &$0.9965$ &  &$38.368$ &$0.00151$ &$0.2033$ &$0.9593$  & &$41.2524$ &$0.00078$ &$0.2601$ &$0.9488$  \\
		&  &$500/50$  &$27.811$  &$0.00016$  &$0.1916$  &$0.9979$  &  &$29.579$ &$0.00179$ &$0.1678$ &$0.9603$ & &$29.0870$ &$0.00090$ &$0.1995$ &$0.9528$  \\
		&  &$500/100$  &$57.927$  &$0.00159$  &$0.4043$  &$0.9979$  &  &$56.712$ &$0.00073$ &$0.4300$ &$0.9610$ & &$56.3514$ &$0.00128$ &$0.3583$ &$0.9576$   \\
		\hline
		&  &\textbf{Average}  &\textbf{37.791}  &\textbf{0.00134}  &\textbf{0.2049}  &\textbf{0.9972}  &  &\textbf{35.953} &\textbf{0.00120} &\textbf{0.2291} &\textbf{0.9586} & &\textbf{37.260} &\textbf{0.00127} &\textbf{0.2194} &\textbf{0.9515}  \\
		\hline\hline
		&  &$300/50$  &$13.567$  &$0.00373$  &$0.0862$  &$0.9962$  &  &$9.444$ &$0.00076$ &$0.1005$ &$0.9470$ & &$11.8334$ &$0.00081$ &$0.1091$ &$0.9344$   \\
		&$0.6$  &$300/100$  &$19.375$  &$0.00106$  &$0.1602$  &$0.9963$  &  &$24.410$ &$0.00101$ &$0.2529$ &$0.9431$  & &$21.4065$ &$0.00114$ &$0.2891$ &$0.9286$   \\
		&  &$500/50$  &$18.388$  &$0.00109$  &$0.3731$  &$0.9977$  &  &$19.895$ &$0.00128$ &$0.2111$ &$0.9207$ & &$18.7170$ &$0.00131$ &$0.2399$ &$0.9411$  \\
		&  &$500/100$  &$36.946$  &$0.00170$  &$0.4550$  &$0.9976$  &  &$38.435$ &$0.00091$ &$0.7091$ &$0.9389$ & &$35.9315$ &$0.00082$ &$0.4335$ &$0.9463$  \\
		\hline
		&  &\textbf{Average}  &\textbf{22.069}  &\textbf{0.00189}  &\textbf{0.2686}  &\textbf{0.9970}  & &\textbf{23.042} &\textbf{0.00099} &\textbf{0.3184} &\textbf{0.9374}  & &\textbf{21.972} &\textbf{0.00102} &\textbf{0.2679} &\textbf{0.9376}   \\
		\hline\hline
		&  &$300/50$  &$1.047$  &$0.00131$  &$0.3326$  &$0.9957$  &  &$-0.334$ &$0.00054$ &$0.1514$ &$0.9048$  & &$1.8247$ &$0.00086$ &$0.1958$ &$0.8991$ \\
		&$0.4$  &$300/100$  &$2.023$  &$0.00166$  &$0.3753$  &$0.9957$  &  &$3.368$ &$0.00085$ &$0.4862$ &$0.8898$  & &$1.6393$ &$0.00100$ &$0.2977$ &$0.8980$ \\
		&  &$500/50$  &$7.989$  &$0.00091$  &$0.7267$  &$0.9973$  &  &$8.896$ &$0.00098$ &$0.2488$ &$0.9136$  & &$9.7090$ &$0.00092$ &$0.3167$ &$0.9069$ \\
		&  &$500/100$  &$16.735$  &$0.00076$  &$1.9874$  &$0.9971$  &  &$17.327$ &$0.00081$ &$0.5745$ &$0.9196$ & &$16.4895$ &$0.00078$ &$0.8149$ &$0.9100$  \\
		\hline
		&   &\textbf{Average} &\textbf{6.948}  &\textbf{0.00116}  &\textbf{0.8555}  &\textbf{0.9964}  &  &\textbf{7.314} &\textbf{0.00080} &\textbf{0.3652} &\textbf{0.9070} & &\textbf{7.416} &\textbf{0.00089} &\textbf{0.4063} &\textbf{0.9035}    \\
		\hline\hline
		&  &$300/50$  &$-9.520$  &$0.00067$  &$1.7936$  &$0.9919$  &  &$-10.415$ &$0.00070$ &$0.2977$ &$0.7858$ & &$-7.2989$ &$0.00074$ &$0.2253$ &$0.7953$  \\
		&$0.2$  &$300/100$  &$-19.515$  &$0.00097$  &$1.4411$  &$0.9893$  &  &$-17.910$ &$0.00066$ &$1.2367$ &$0.7744$  & &$-17.1456$ &$0.00066$ &$1.2029$ &$0.7786$ \\
		&  &$500/50$  &$-1.981$  &$0.00075$  &$1.8656$  &$0.9957$  &  &$-0.521$ &$0.00082$ &$0.5231$ &$0.8103$ & &$-0.2167$ &$0.00072$ &$0.6751$ &$0.8152$  \\
		&  &$500/100$  &$-3.107$  &$0.00097$  &$3.3941$  &$0.9942$  &  &$-0.292$ &$0.00070$ &$1.2026$ &$0.8238$ & &$0.0385$ &$0.00073$ &$1.4124$ &$0.8006$  \\
		\hline
		&  &\textbf{Average}  &\textbf{-8.531}  &\textbf{0.00084} &\textbf{2.1236}  &\textbf{0.9928}  &  &\textbf{-7.284} &\textbf{0.00072} &\textbf{0.8150} &\textbf{0.7986} & &\textbf{-6.156} &\textbf{0.00071} &\textbf{0.8789} &\textbf{0.7974}    \\
		\hline\hline
	\end{tabular}
	}
	\label{tab: gausRes}
\end{table*}

\begin{table}[h]
	\caption{Results for the proposed MQPAM. Note that the values for $\mu$ are small for MQPAM.}
	\centering
    \resizebox{\columnwidth}{!}{
	\begin{tabular}{ccccccc}
		\hline
		&$\mu$ &$n/p$ &Obj &Err &CPU &Spar \\
		\hline
		& &$300/50$ &$51.0995$ &$0.0023$ &$0.0402$ &$0.9765$ \\
		&$1\times 10^{-10}$ &$300/100$ & $99.5474$ &$0.0019$ &$0.0984$ &$0.9766$  \\
		& &$500/50$ &$49.8216$ &$0.0008$ &$0.1346$ &$0.9780$ \\
		& &$500/100$ &$99.0804$ &$0.0014$ &$0.1819$ &$0.9779$ \\
		\hline
		& &\textbf{Average}  &\textbf{74.887} &\textbf{0.0016} &\textbf{0.1138} &\textbf{0.9773} \\
		\hline\hline
		& &$300/50$ &$49.7482$ &$0.0018$ &$0.0434$ &$0.9766$ \\
		&$1\times 10^{-8}$ &$300/100$ &$99.7909$ &$0.0011$ &$0.1261$ &$0.9766$ \\
		& &$500/50$ &$49.9754$ &$0.0014$ &$0.0882$ &$0.9780$ \\
		& &$500/100$ &$99.1763$ &$0.0006$ &$0.2088$ &$0.9780$ \\
		\hline
		&  &\textbf{Average}  &\textbf{74.673} &\textbf{0.0012} &\textbf{0.1167} &\textbf{0.9773} \\
		\hline\hline
		& &$300/50$ &$49.5219$ &$0.0002$ &$0.0331$ &$0.9766$ \\
		&$1\times 10^{-6}$ &$300/100$ &$99.5395$ &$0.0019$ &$0.0863$ &$0.9765$ \\
		& &$500/50$ &$50.1476$ &$0.0020$ &$0.0816$ &$0.9779$ \\
		& &$500/100$ &$100.5333$  &$0.0019$ &$0.1765$ &$0.9780$ \\
		\hline
		& &\textbf{Average}  &\textbf{74.936} &\textbf{0.0015} &\textbf{0.094} &\textbf{0.9772}  \\
		\hline\hline
		& &$300/50$ &$49.5934$ &$0.0011$ &$0.0541$ &$0.3424$ \\
		&$1\times 10^{-4}$ &$300/100$ &$99.0451$ &$0.0008$ &$0.1295$ &$0.4659$ \\
		& &$500/50$ &$50.5851$ &$0.0023$ &$0.0860$ &$0.3647$ \\
		& &$500/100$ &$99.3861$  &$0.0017$ &$0.1596$ &$0.4283$ \\
		\hline
		& &\textbf{Average}  &\textbf{74.652} &\textbf{0.0015} &\textbf{0.1073} &\textbf{0.4003}  \\
		\hline\hline
		& &$300/50$ &$50.830$ &$0.0014$ &$0.0667$ &$0.1611$ \\
		&$1\times 10^{-2}$ &$300/100$ &$101.111$ &$0.002$ &$0.1304$ &$0.3250$ \\
		& &$500/50$ &$50.653$ &$0.0010$ &$0.0903$ &$0.0991$ \\
		& &$500/100$ &$101.001$  &$0.0005$ &$0.2076$ &$0.1960$ \\
		\hline
		& &\textbf{Average}  &\textbf{75.899} &\textbf{0.0012} &\textbf{0.1238} &\textbf{0.1953}  \\
		\hline\hline
	\end{tabular}
 }
	\label{tab:resMam}
\end{table}

\section{Sparse Principal Component Analysis}

We proceed by applying the MQPAM to solve the SPCA problem 
\begin{equation}
\begin{aligned}
&\underset{\mathbf{X}}{minimize}\, F\left( \mathbf{X} \right)= - \frac{\mu}{2}\text{Tr}\left( \mathbf{X}^\top \mathbf{A}^\top \mathbf{A}\mathbf{X} \right) + \| \mathbf{X} \|_1 \\
&subject\,to\quad \mathbf{X}\in \text{St}\left(n,\,p \right),
\end{aligned}
\label{eq:spca}
\end{equation}
which is an instance of (\ref{eq:spca1}) with a slight difference that the regularization parameter $\mu$ is coupled on the first term as to be consistent with the formulations in \cite{wang2008new,xie2018new}. As we shall see in the experiments, the value of $\mu$ from this slight difference still acts as a regularization parameter which imposes sparsity to the solutions.

To apply MQPAM, reformulate problem (\ref{eq:spca}) into (\ref{eq: probSplit1}) by adding a quadratic penalty,
\begin{equation}
\begin{aligned}
&\underset{\mathbf{X}}{minimize}\, - \frac{\mu}{2}\text{Tr}\left(\mathbf{X}^\top \mathbf{A}^\top \mathbf{AX} \right) + \| \mathbf{Y} \|_1 + \frac{\beta}{2}\|\mathbf{Y} - \mathbf{X} \|_2^2 \\
&\text{\textit{subject to}}\quad \mathbf{Y}= \mathbf{X},\, \mathbf{X}\in \text{St}\left(n,\,p \right). \\
\end{aligned}
\label{eq:probRqpam_PCA}
\end{equation}
RQPAM then minimizes $\mathbf{X}$ and $\mathbf{Y}$ alternately as,
\begin{equation}
\begin{cases}
\mathbf{X}_{k+1} = \underset{\mathbf{X} \in \text{St}\left(n,\,p \right)}{arg\,min}\, - \frac{\mu}{2}\text{Tr}\left(\mathbf{X}^\top \mathbf{A}^\top \mathbf{AX} \right) + \frac{\beta}{2}\| \mathbf{Y}_k - \mathbf{X} \|_2^2 \\
\mathbf{Y}_{k+1} = \underset{\mathbf{Y}}{arg\,min}\, \frac{\beta}{2}\| \mathbf{Y} - \mathbf{X}_{k+1} \|_2^2 + \| \mathbf{Y}\|_1.
\end{cases}
\end{equation}
The manifold-constrained subproblem $\mathbf{X}$, can be solved using several iterations of the RGD as in (\ref{eq:rgd}) with Euclidean gradient 
$\mathbf{G}_k = - \mu \mathbf{A}^\top \mathbf{AX} + \beta \left( \mathbf{X} - \mathbf{Y} \right)$.

The $\mathbf{Y}$ subproblem is an $\ell_1$-norm minimization problem in the form of the proximal operator (\ref{eq:prox_lam}) thus, having a closed-form solution via the soft-thresholding as follows
\begin{equation}
\text{soft}\left( \mathbf{X}_{k+1},\,\frac{1}{\beta}\right) = \text{sgn}\left( \mathbf{X}_{k+1}\right) \odot \text{max} \left( \left| \mathbf{X}_{k+1}\right| - \frac{1}{\beta},0 \right).
\end{equation}
The overall algorithm of MQPAM applied to SPCA is listed as Algorithm \ref{al:rqpam_pca}. 
\begin{algorithm}
\caption{MQPAM for Sparse PCA problem (\ref{eq:probRqpam_PCA}) }
\label{al:rqpam_pca}

\KwInit{$\mathbf{X}_0$, $\mathbf{Y}_0$, $\eta > 0$, $tol = 1\times10^{-8}$, $Nit$, $Rgdit$}

\For{$k=0$ \KwTo $Nit$}
{
    \For{$j=0$ \KwTo $Rgdit$}
    {
        $\mathbf{G}_{k+1} = - \mu \mathbf{A}^\top \mathbf{AX}_k + \beta \left( \mathbf{X}_k - \mathbf{Y}_k \right)$ \\
        $\mathcal{G}_{k+1} = \text{Proj}_{T_X \mathcal{M}}\left(\mathbf{G}_k\right)$ \\
        \If{$\| \mathcal{G}_{k+1} \|_2 < tol$}
        {
            Break
        }
        $\mathbf{X}_{k+1} = \text{Retr}_X \left( -\eta \mathcal{G}_{k+1} \right)$ \\
    }

    $\mathbf{Y}_{k+1} =  \text{sgn}\left( \mathbf{X}_{k+1}\right) \odot \text{max} \left( \left| \mathbf{X}_{k+1}\right| - \frac{1}{\beta},0 \right)$

$k= k+1$
}
\end{algorithm}

\section{Experimental Results}
\label{sec:ER}
We conducted experiments\footnote{Conducted with MATLAB\textsuperscript{\textregistered} 2023 on 64-bit Intel \textsuperscript{\textregistered} Core\texttrademark  i5 2.3GHz CPU with 8GB RAM. Code will be released at: \href{https://github.com/tarmiziAdam2005}{https://github.com/tarmiziAdam2005} } using MQPAM to the SPCA problem (\ref{eq:spca}) and compared the MQPAM with several other methods. The compared methods are Splitting with Orthogonality Constraint (SOC) \cite{lai2014splitting}, MADMM \cite{kovnatsky2016madmm}, and the recently proposed RADMM \cite{li2022riemannian}. All the compared methods are operator-splitting methods that were used in the SPCA problem. Regarding the parameters, for the SOC, we set the proximal gradient stepsize $\eta=1 \times 10^{-2}$ and $\rho = 50$. For MADMM, $\eta= 1 \times 10^{-2}$ for the RGD step and $\rho = 100$. For RADMM, we set $\eta = 1 \times 10^{-2}$, $\rho = 100$, and $\gamma = 1 \times 10^{-8}$. For more details on these parameters, the reader is referred to their respective papers. For MQPAM, we set $\eta = 1 \times 10^{-2}$ and $\beta = 100.5$. The RGD iterations were set to 100 iterations for MADMM, RADMM, and MQPAM We also used the relative error-stopping criterion for all the algorithms as follows:
\[
\frac{\| \mathbf{X}_{k+1} - \mathbf{X}_k\|}{\| \mathbf{X}_{k+1} \|_2} \leq 1 \times 10^{-5}.
\]

The data matrix $\mathbf{A}\in \mathbb{R}^{m \times n}$ was randomly generated having entries following a Gaussian distribution. For the Stiefel manifold $St\left(n,p \right)$, we experimented with values of $n = 300$ and $n = 500$ while $p = 50$ and $p=100$. Finally, all the experiments for each case are repeated 50 times and the average is reported.

Tables \ref{tab: gausRes} and \ref{tab:resMam} show the comparative results for all the algorithms. In the tables, we report the objective function value "Obj" ( $F\left( \mathbf{X}_k  \right)$), relative error "Err", CPU time (seconds), and the sparsity level which is the percentage of the zero entries of the solution for different values of $\mu$.

From Tables \ref{tab: gausRes} and \ref{tab:resMam}, it is observed that MQPAM does not give a small objective value compared to the other methods however, the relative error value is comparable to other methods.

Overall, MQPAM is faster compared to the other methods while giving higher sparsity levels compared to MADMM and RADMM (for $\mu = 1\times 10^{-6}$ to $1 \times 10^{-10}$). This shows that the simple steps of the MQPAM reduce the per-iteration complexity and hence, the CPU time.
On another note, the fast computation time is probably due to no dual update step in the MQPAM iteration, unlike SOC, MADMM, and RADMM.
However, it is observed that the SOC method gives a higher sparsity level for all the compared methods.

In Table \ref{tab:resMam} we can observe that the MQPAM is very sensitive to the sparsity parameter $\mu$, unlike the other methods. This could be seen in the sudden drop of sparsity level when $\mu = 1\times 10^{-4}$ and $\mu = 1\times 10^{-2}$. Care should be taken in choosing this parameter to obtain reasonable results for the MQPAM. 

\section{Conclusion}
\label{sec:C}
In this paper, we proposed a manifold quadratic penalty alternating minimization (MQPAM) algorithm for optimization on the Stiefel manifold. The MQPAM is very simple to implement. The subproblems of the MQPAM only require an RGD step and a closed-form proximal step. Experimental results for SPCA show that the MQPAM is fast in terms of CPU time while achieving comparable or better sparsity results for the solution. As a future endeavor, we are currently studying the detailed convergence analysis of the MQPAM and other possible applications of the MQPAM.

\bibliographystyle{IEEEtran}
\bibliography{EUSIPCO24}

\begin{thebibliography}{10}
\providecommand{\url}[1]{#1}
\csname url@samestyle\endcsname
\providecommand{\newblock}{\relax}
\providecommand{\bibinfo}[2]{#2}
\providecommand{\BIBentrySTDinterwordspacing}{\spaceskip=0pt\relax}
\providecommand{\BIBentryALTinterwordstretchfactor}{4}
\providecommand{\BIBentryALTinterwordspacing}{\spaceskip=\fontdimen2\font plus
\BIBentryALTinterwordstretchfactor\fontdimen3\font minus
  \fontdimen4\font\relax}
\providecommand{\BIBforeignlanguage}[2]{{%
\expandafter\ifx\csname l@#1\endcsname\relax
\typeout{** WARNING: IEEEtran.bst: No hyphenation pattern has been}%
\typeout{** loaded for the language `#1'. Using the pattern for}%
\typeout{** the default language instead.}%
\else
\language=\csname l@#1\endcsname
\fi
#2}}
\providecommand{\BIBdecl}{\relax}
\BIBdecl

\bibitem{boumal2023introduction}
N.~Boumal, \emph{An introduction to optimization on smooth manifolds}.\hskip
  1em plus 0.5em minus 0.4em\relax Cambridge University Press, 2023.

\bibitem{absil2008optimization}
P.-A. Absil, R.~Mahony, and R.~Sepulchre, \emph{Optimization algorithms on
  matrix manifolds}.\hskip 1em plus 0.5em minus 0.4em\relax Princeton
  University Press, 2008.

\bibitem{yang2011}
Y.~Yang, H.~T. Shen, Z.~Ma, Z.~Huang, and X.~Zhou, ``$\ell_{2,1}$-norm
  regularized discriminative feature selection for unsupervised learning,'' in
  \emph{IJCAI international joint conference on artificial intelligence}, 2011.

\bibitem{tang2012}
J.~Tang and H.~Liu, ``Unsupervised feature selection for linked social media
  data,'' in \emph{Proceedings of the 18th ACM SIGKDD international conference
  on Knowledge discovery and data mining}, 2012, pp. 904--912.

\bibitem{sun2016b}
J.~Sun, Q.~Qu, and J.~Wright, ``Complete dictionary recovery over the sphere i:
  Overview and the geometric picture,'' \emph{IEEE Transactions on Information
  Theory}, vol.~63, no.~2, pp. 853--884, 2016.

\bibitem{sun2016a}
------, ``Complete dictionary recovery over the sphere ii: Recovery by
  riemannian trust-region method,'' \emph{IEEE Transactions on Information
  Theory}, vol.~63, no.~2, pp. 885--914, 2016.

\bibitem{zhang2017global}
Y.~Zhang, Y.~Lau, H.-w. Kuo, S.~Cheung, A.~Pasupathy, and J.~Wright, ``On the
  global geometry of sphere-constrained sparse blind deconvolution,'' in
  \emph{Proceedings of the IEEE Conference on Computer Vision and Pattern
  Recognition}, 2017, pp. 4894--4902.

\bibitem{ozolicnvs2013}
V.~Ozoli{\c{n}}{\v{s}}, R.~Lai, R.~Caflisch, and S.~Osher, ``Compressed modes
  for variational problems in mathematics and physics,'' \emph{Proceedings of
  the National Academy of Sciences}, vol. 110, no.~46, pp. 18\,368--18\,373,
  2013.

\bibitem{jolliffe2003}
I.~T. Jolliffe, N.~T. Trendafilov, and M.~Uddin, ``A modified principal
  component technique based on the lasso,'' \emph{Journal of computational and
  Graphical Statistics}, vol.~12, no.~3, pp. 531--547, 2003.

\bibitem{huang2022}
W.~Huang and K.~Wei, ``An extension of fast iterative shrinkage-thresholding
  algorithm to riemannian optimization for sparse principal component
  analysis,'' \emph{Numerical Linear Algebra with Applications}, vol.~29,
  no.~1, p. e2409, 2022.

\bibitem{xiao2021exact}
N.~Xiao, X.~Liu, and Y.-x. Yuan, ``Exact penalty function for $\ell_{2,1}$-norm
  minimization over the stiefel manifold,'' \emph{SIAM Journal on
  Optimization}, vol.~31, no.~4, pp. 3097--3126, 2021.

\bibitem{pearson1901}
K.~Pearson, ``Liii. on lines and planes of closest fit to systems of points in
  space,'' \emph{The London, Edinburgh, and Dublin philosophical magazine and
  journal of science}, vol.~2, no.~11, pp. 559--572, 1901.

\bibitem{hotelling1933}
H.~Hotelling, ``Analysis of a complex of statistical variables into principal
  components.'' \emph{Journal of educational psychology}, vol.~24, no.~6, p.
  417, 1933.

\bibitem{hu2020brief}
J.~Hu, X.~Liu, Z.-W. Wen, and Y.-X. Yuan, ``A brief introduction to manifold
  optimization,'' \emph{Journal of the Operations Research Society of China},
  vol.~8, pp. 199--248, 2020.

\bibitem{sato2021riemannian}
H.~Sato, \emph{Riemannian optimization and its applications}.\hskip 1em plus
  0.5em minus 0.4em\relax Springer, 2021, vol. 670.

\bibitem{huang2022riemannian}
W.~Huang and K.~Wei, ``Riemannian proximal gradient methods,''
  \emph{Mathematical Programming}, vol. 194, no. 1-2, pp. 371--413, 2022.

\bibitem{chen2020proximal}
S.~Chen, S.~Ma, A.~Man-Cho~So, and T.~Zhang, ``Proximal gradient method for
  nonsmooth optimization over the stiefel manifold,'' \emph{SIAM Journal on
  Optimization}, vol.~30, no.~1, pp. 210--239, 2020.

\bibitem{huang2021robust}
M.~Huang, S.~Ma, and L.~Lai, ``Robust low-rank matrix completion via an
  alternating manifold proximal gradient continuation method,'' \emph{IEEE
  Transactions on Signal Processing}, vol.~69, pp. 2639--2652, 2021.

\bibitem{kovnatsky2016madmm}
A.~Kovnatsky, K.~Glashoff, and M.~M. Bronstein, ``Madmm: a generic algorithm
  for non-smooth optimization on manifolds,'' in \emph{Computer Vision--ECCV
  2016: 14th European Conference, Amsterdam, The Netherlands, October 11-14,
  2016, Proceedings, Part V 14}.\hskip 1em plus 0.5em minus 0.4em\relax
  Springer, 2016, pp. 680--696.

\bibitem{li2022riemannian}
J.~Li, S.~Ma, and T.~Srivastava, ``A riemannian admm,'' \emph{arXiv preprint
  arXiv:2211.02163}, 2022.

\bibitem{wang2008new}
Y.~Wang, J.~Yang, W.~Yin, and Y.~Zhang, ``A new alternating minimization
  algorithm for total variation image reconstruction,'' \emph{SIAM Journal on
  Imaging Sciences}, vol.~1, no.~3, pp. 248--272, 2008.

\bibitem{xie2018new}
J.~Xie, A.~Liao, and Y.~Lei, ``A new accelerated alternating minimization
  method for analysis sparse recovery,'' \emph{Signal Processing}, vol. 145,
  pp. 167--174, 2018.

\bibitem{xu2011image}
L.~Xu, C.~Lu, Y.~Xu, and J.~Jia, ``Image smoothing via l 0 gradient
  minimization,'' in \emph{Proceedings of the 2011 SIGGRAPH Asia conference},
  2011, pp. 1--12.

\bibitem{lai2014splitting}
R.~Lai and S.~Osher, ``A splitting method for orthogonality constrained
  problems,'' \emph{Journal of Scientific Computing}, vol.~58, pp. 431--449,
  2014.

\bibitem{siegel2021accelerated}
J.~W. Siegel, ``Accelerated optimization with orthogonality constraints,''
  \emph{Journal of Computational Mathematics}, vol.~39, no.~2, p. 207, 2021.

\bibitem{chen2016augmented}
W.~Chen, H.~Ji, and Y.~You, ``An augmented lagrangian method for
  $\ell_1$-regularized optimization problems with orthogonality constraints,''
  \emph{SIAM Journal on Scientific Computing}, vol.~38, no.~4, pp. B570--B592,
  2016.

\bibitem{boyd2011distributed}
S.~Boyd, N.~Parikh, E.~Chu, B.~Peleato, J.~Eckstein \emph{et~al.},
  ``Distributed optimization and statistical learning via the alternating
  direction method of multipliers,'' \emph{Foundations and
  Trends{\textregistered} in Machine learning}, vol.~3, no.~1, pp. 1--122,
  2011.

\bibitem{goldstein2009split}
T.~Goldstein and S.~Osher, ``The split bregman method for l1-regularized
  problems,'' \emph{SIAM journal on imaging sciences}, vol.~2, no.~2, pp.
  323--343, 2009.

\bibitem{beck2009fast}
A.~Beck and M.~Teboulle, ``A fast iterative shrinkage-thresholding algorithm
  for linear inverse problems,'' \emph{SIAM journal on imaging sciences},
  vol.~2, no.~1, pp. 183--202, 2009.

\bibitem{deng2023manifold}
K.~Deng and Z.~Peng, ``A manifold inexact augmented lagrangian method for
  nonsmooth optimization on riemannian submanifolds in euclidean space,''
  \emph{IMA Journal of Numerical Analysis}, vol.~43, no.~3, pp. 1653--1684,
  2023.

\bibitem{tan2014smoothing}
Z.~Tan, Y.~C. Eldar, A.~Beck, and A.~Nehorai, ``Smoothing and decomposition for
  analysis sparse recovery,'' \emph{IEEE Transactions on Signal Processing},
  vol.~62, no.~7, pp. 1762--1774, 2014.

\bibitem{courant1943}
R.~Courant, ``Variational methods for the solution of problems of equilibrium
  and vibrations,'' 1943.

\bibitem{wright2022optimization}
S.~J. Wright and B.~Recht, \emph{Optimization for data analysis}.\hskip 1em
  plus 0.5em minus 0.4em\relax Cambridge University Press, 2022.

\bibitem{boumal2019}
N.~Boumal, P.-A. Absil, and C.~Cartis, ``Global rates of convergence for
  nonconvex optimization on manifolds,'' \emph{IMA Journal of Numerical
  Analysis}, vol.~39, no.~1, pp. 1--33, 2019.

\bibitem{bento17}
G.~C. Bento, O.~P. Ferreira, and J.~G. Melo, ``Iteration-complexity of
  gradient, subgradient and proximal point methods on riemannian manifolds,''
  \emph{Journal of Optimization Theory and Applications}, vol. 173, no.~2, pp.
  548--562, 2017.

\end{thebibliography}

\end{document}